\documentclass{article}

\usepackage{arxiv}

\usepackage[utf8]{inputenc} 
\usepackage[T1]{fontenc}    
\usepackage{hyperref}       
\usepackage{url}            
\usepackage{booktabs}       
\usepackage{amsfonts}       
\usepackage{nicefrac}       
\usepackage{microtype}      
\usepackage{lipsum}

\usepackage{mathptmx,graphicx}

\title{Novel Numerical Algorithm with Fourth-Order Accuracy for the Direct Zakharov-Shabat Problem}

\author{
  Sergey Medvedev$^{1,2,*}$, Irina Vaseva$^{1}$, Igor Chekhovskoy$^{1,2}$, Mikhail Fedoruk$^{1,2}$\\
$^{1}$ Institute of Computational Technologies, SB RAS, Novosibirsk 630090, Russia,\\
$^{2}$ Novosibirsk State University, Novosibirsk 630090, Russia,\\
* Corresponding author: medvedev@ict.nsc.ru
}

\begin{document}
\maketitle

\begin{abstract}
We propose a new high-precision algorithm for solving the initial problem for the Zakharov-Shabat system. This method has the fourth order of accuracy and is a generalization of the second order Boffetta-Osborne scheme. It is allowed by our method to solve more effectively the Zakharov-Shabat spectral problem for continuous and discrete spectra.
\end{abstract}


The solution of the direct problem for the Zakharov-Shabat problem (ZSP) is the first step in the inverse scattering transform (IST) for solving the nonlinear Schr\"{o}dinger equation (NLSE)~\cite{ZakharovShabat1972}.
The numerical implementation of the IST has gained great importance and attracted special attention since Hasegawa and Tappert~\cite{Hasegawa1973a} proposed to use soliton solutions as a bit of information for fiber optic data transmission.

The direct scattering problem is solved by spectral data. To calculate them it is necessary to solve the initial problem for the Zakharov-Shabat system. Therefore, a lot of effort was made to find effective numerical methods for solving this problem. An overview of the methods used can be found in~\cite{Yousefi2014II,Turitsyn2017Optica,Vasylchenkova2017a}. Currently, one of the most effective methods for solving the ZSP is the Boffetta-Osborne (BO) method~\cite{Boffetta1992a}, which has the second order of approximation. Comparisons of this method with other methods were carried out in~\cite{Vasylchenkova2017a,Burtsev1998}.

Besides the approximation accuracy, it is necessary to have an algorithm  requiring a minimum computational time to get a discrete set of spectral parameters with sufficient accuracy.
This direction is implemented in the fast algorithm (FNFT) for solving the direct ZSP using the modified Ablowitz-Ladik method~\cite{Wahls2013, Turitsyn2017Optica}. The BO method does not allow a direct application of the fast algorithm. But the fast method can be applied to the exponential approximation of the transition matrix of the BO method~\cite{Prins2018a}. In this Letter, we will focus on building the method of the fourth order of accuracy on a uniform grid. For a non-uniform grid, a fourth order scheme~\cite{Blanes2017} was applied in~\cite{Prins2018}. In perspective, the exponential approximation can be applied to our scheme so that we can use the fast algorithm.


We write the Zakharov-Shabat system in a matrix form
\begin{equation}\label{psit}
\frac{d}{dt}{ \Psi}(t)=Q(t){\Psi}(t),
\end{equation}
where $\Psi(t)$ is a complex vector function of the real argument~$t$, $Q(t)$ is a complex matrix
$$
{\Psi}(t) = \left(
\begin{array}{c}
    \psi_1(t)\\\psi_2(t)
\end{array} \right),\quad
Q(t) = \left(
\begin{array}{cc}
    -i\zeta&q(t,z_0)\\-\sigma q^*(t,z_0)&i\zeta
\end{array} \right),
$$
where $\sigma=\pm 1$ for anomalous and normal dispersion, $z_0$ plays the role of a parameter and will not be used further. The asterisk means the complex conjugation.

Consider the Jost initial conditions
\begin{equation}\label{psi0}
\left(
\begin{array}{c}
    \psi_{1}\\\psi_{2}
\end{array}
\right) = \left(
\begin{array}{c}
    e^{-i\zeta t}\\0
\end{array}
\right)[1+o(1)],\quad t\to-\infty,
\end{equation}
which define the Jost solutions for real~$\zeta=\xi$.
The coefficients of the scattering matrix~$a(\xi)$ and~$b(\xi)$ are obtained as limits
\begin{equation}\label{ab}
a(\xi)=\lim_{t\to\infty}\,\psi_1(t,\xi)\,e^{i\xi t},\quad b(\xi)=\lim_{t\to\infty}\,\psi_2(t,\xi)\,e^{-i\xi t}.
\end{equation}
The function~$a(\xi)$ can be extended to the upper half-plane $\xi\to \zeta$, where~$\zeta$ is a complex number with the positive imaginary part $\eta=\mbox{Im}\,\zeta>0$.
The spectral data are determined by $a(\zeta)$ and $b(\zeta)$ in the following way:\\
(1) the zeros of $a(\zeta)=0$ define the discrete spectrum $\{\zeta_k\}$, $k=1,...,K$ of ZSP~(\ref{psit}) and phase coefficients
$$r_k=\left.\frac{b(\zeta)}{a'(\zeta)}\right|_{\zeta=\zeta_k},\quad\mbox{where}\quad a'(\zeta)=\frac{da(\zeta)}{d\zeta};$$
(2) the continuous spectrum is determined by the reflection coefficient
$$r(\xi)=\frac{b(\xi)}{a(\xi)},\quad \xi\in\mathbb{R}. $$

The matrix~$Q(t)$ in the system~(\ref{psit}) becomes the skew-Her\-mi\-tian ($Q^*=-Q^T$) when the spectral parameter $\zeta=\xi$ is real and $\sigma=1$. Therefore, the system~(\ref{psit}) preserves the integral
\begin{equation}
\frac{d}{dt}\left(|\psi_1(t)|^2+|\psi_2(t)|^2\right)=0.
\end{equation}
Taking into account the boundary conditions~(\ref{psi0}), we have
\begin{equation}\label{conser}
|\psi_1(t)|^2+|\psi_2(t)|^2=1.
\end{equation}

In addition, the trace formula is valid~\cite{Ablowitz1981}
\begin{eqnarray}
C_n=-\frac{1}{\pi}\int\limits_{-\infty}^\infty\,(2i\xi)^n\,\ln|a(\xi)|^2\,d\xi+
+\sum\limits_{k=1}^K\,\frac{1}{(n+1)}
\left[(2i\zeta_k^*)^{n+1}-(2i\zeta_k)^{n+1}\right],\nonumber
\end{eqnarray}
which connects the NLSE integrals $C_n$ with the coefficient $a(\xi)$ and the discrete spectrum $\zeta_k$. The first integrals have the form
$$C_0=\int\limits_{-\infty}^\infty|q|^2dt,\enskip C_1=\int\limits_{-\infty}^\infty qq^*_tdt, \enskip C_2=\int\limits_{-\infty}^\infty(qq^*_{tt}+|q|^4)dt,$$
$$C_3=\int\limits_{-\infty}^\infty(qq^*_{ttt}+4|q|^2qq^*_t+|q|^2q^*q_t)dt. $$
This formula with $n=0$ is called the Parseval nonlinear equality and is used to verify the numerical calculations and the consistency of the continuous and discrete spectra found.


We solve the system~(\ref{psit}). The matrix~$Q(t)$ linearly depends on the complex function~$q(t)$, which is given in the whole nodes of the uniform grid $t_n=-L+\tau n$ with a step~$\tau$ on the interval $[-L,L]$. If the total number of points is~$2M+1$, then the grid step is $\tau=L/M$.
Since the matrix~$Q(t)$ is specified only on the grid, Boffetta and Osborne suggested replacing the original system on the interval $[t_n-\frac{\tau}{2},t_n+\frac{\tau}{2}]$ with an approximate system with constant coefficients \cite{Boffetta1992a}
\begin{equation}
\frac{d}{dt}\Psi(t)=Q(t_n)\Psi(t),\quad Q_n=Q(t_n),
\end{equation}
which is easily solved on the selected interval and gives the transition matrix from the layer $n-\frac{1}{2}$ to the layer $n+\frac{1}{2}$:
\begin{equation}\label{T0}
\Psi_{n+\frac{1}{2}}=e^{\tau Q_n}\Psi_{n-\frac{1}{2}}.
\end{equation}
This method has proven itself well, but nonetheless one would like to get a more accurate solution.
So we formulate our task: it is required on the interval $[t_n-\frac{\tau}{2},t_n+\frac{\tau}{2}]$ to build the transition matrix from $\Psi_{n-\frac{1}{2}}$ to $\Psi_{n+\frac{1}{2}}$ with maximum accuracy and minimum computational cost.

The first step is a change of variables
$$\Psi(t)=e^{tQ_n}Y(t),$$
so the initial system takes the form
\begin{equation}\label{yt}
\frac{d}{dt}{Y}(t)=L(t)Y(t),\quad L(t)=e^{-tQ_n}(Q(t)-Q_n)e^{tQ_n}.
\end{equation}
In this form, the linear matrix~$L(t)$ becomes zero at~$t=t_n$, and the derivative of~$Y(t)$ is zero at this point. This means that the solution is almost constant in the neighborhood of~$t_n$. If the original system is replaced by an approximate
\begin{equation}\label{app0}
\frac{d}{dt}{Y}(t)=0,
\end{equation}
then the transition from $Y_{n-\frac{1}{2}}$ to $Y_{n+\frac{1}{2}}$ becomes trivial:
\begin{equation}\label{YY}
Y_{n+\frac{1}{2}}=Y_{n-\frac{1}{2}}.
\end{equation}
Returning to the original values of the variable~$\Psi$ at the points $t_n-\frac{\tau}{2}$ and $t_n+\frac{\tau}{2}$, we get
\begin{equation}\label{psiy}
 \Psi_{n-\frac{1}{2}}=e^{\left(t_n-\frac{\tau}{2}\right)Q_n}Y_{n-\frac{1}{2}},\quad \Psi_{n+\frac{1}{2}}=e^{\left(t_n+\frac{\tau}{2}\right)Q_n}Y_{n+\frac{1}{2}},
\end{equation}
which, taking into account (\ref{YY}), exactly gives a transition in the BO scheme (\ref{T0}). Since we are interested in the values of~$Y$ only in the grid nodes, the solution~(\ref{YY}) can be interpreted as a solution to the difference equation
\begin{equation}\label{R0}
\frac{Y_{n+\frac{1}{2}}-Y_{n-\frac{1}{2}}}{\tau}=L_n\frac{Y_{n+\frac{1}{2}}+Y_{n-\frac{1}{2}}}{2},
\end{equation}
which is an approximation of the continuous equation~(\ref{yt}) given that $L_n=L(t_n)=0$. By decomposing~(\ref{R0}) into a Taylor series at the point~$t=t_n$, we get the second order of approximation
$$\frac{d}{dt}Y\left(t_n\right)-L(t_n)Y(t_n)\approx \frac{\tau^2}{24}\frac{d^3}{dt^3}Y(t_n).$$
Thus, we have shown that the BO scheme corresponds to the simplest finite-difference approximation~(\ref{R0}).


There are two possibilities to construct more complex approximations for the equation~(\ref{yt}) on the interval $[t_n-\frac{\tau}{2}, t_n+\frac{\tau}{2}]$. The first is to build a finite difference analog for this equation. The second possibility is to construct an approximation of the operator~$L(t)$ on the entire interval $[t_n-\frac{\tau}{2}, t_n+\frac{\tau}{2}]$ according to the existing values of~$Q_n$ on a regular grid and the subsequent solution of such a system by any analytical method.

Consider the first approach. Since we want to refine the BO scheme, we take the function~$Y$ only in two nodes of the grid $Y_{n-\frac{1}{2}}$ and $Y_{n+\frac{1}{2}}$. For the matrix~$L$, we take the three nearest values~$L_{n-1}$, $L_n$ and $L_{n+1}$. Using these values, we will look for a scheme using the method of uncertain coefficients
\begin{eqnarray}\label{abcd}
&&\frac{Y_{n+\frac{1}{2}}-Y_{n-\frac{1}{2}}}{\tau}=
\left(\alpha L_{n+1}+\beta L_{n-1}\right)Y_{n+\frac{1}{2}}+\left(\gamma L_{n+1}
+\delta L_{n-1}\right)Y_{n-\frac{1}{2}}.\nonumber
\end{eqnarray}
Here we used the condition~$L_n=0$, to drop the terms with~$L_n$ in the right-hand side. By decomposing~(\ref{abcd}) into a Taylor series at the point~$t=t_n$ and using the equation~(\ref{yt}) at this point and its time derivatives, we get that the expression~(\ref{abcd}) has at least the fourth order approximation in $\tau$:
$$\frac{d}{dt}{Y}(t_n)-L(t_n)Y(t_n)\approx\frac{\tau^4}{24}(\beta -\alpha)\frac{d{L}_n}{dt}\frac{d^2L_n}{dt^2}Y_n-\frac{\tau^4}{5760}\left(17\frac{d^4L_n}{dt^4}+12\frac{d^2L_n}{dt^2}\frac{dL_n}{dt}\right)Y_n$$
for
$$\gamma =\frac{1}{24}-\alpha,\quad \delta=\frac{1}{24}-\beta$$
and arbitrary $\alpha$ and $\beta$. The resulting scheme can be rewritten as
\begin{eqnarray}\label{sch}
    &\left[I-\tau \alpha L_{n+1}-\tau \beta L_{n-1}\right]Y_{n+\frac{1}{2}}=\left[I+\tau\left(\frac{1}{24}-\alpha \right)L_{n+1}+
    \tau\left(\frac{1}{24}-\beta \right)L_{n-1}\right]Y_{n-\frac{1}{2}},\nonumber
\end{eqnarray}
where
$$L_{n+1}=e^{-(t_n+\tau) Q_n}\left(Q_{n+1}-Q_n\right)e^{(t_n+\tau) Q_n},$$
$$ L_{n-1}=e^{-(t_n-\tau) Q_n}\left(Q_{n-1}-Q_n\right)e^{(t_n-\tau) Q_n}.$$
In the original variables, this scheme will take the form
\begin{eqnarray}\label{schpsiM}
 &\left[I-\tau \alpha M_{n+1}-\tau \beta M_{n-1}\right]e^{-\frac{\tau}{2}Q_n}\Psi_{n+\frac{1}{2}}=
 \left[I+\tau\left(\frac{1}{24}-\alpha \right)M_{n+1}+
 \tau\left(\frac{1}{24}-\beta \right)M_{n-1}\right]e^{\frac{\tau}{2}Q_n}\Psi_{n-\frac{1}{2}},\nonumber
\end{eqnarray}
where
$$M_{n+1}=e^{-\tau Q_n}\left(Q_{n+1}-Q_n\right)e^{\tau Q_n},$$
$$ M_{n-1}=e^{\tau Q_n}\left(Q_{n-1}-Q_n\right)e^{-\tau Q_n}.$$

For real values $\zeta=\xi$, energy conservation~(\ref{conser}) is important, so if $\alpha  = \beta = 1/48$, then the transition operator
\begin{eqnarray}\label{schpsiMS}
T=e^{\frac{\tau}{2}Q_n}\left[I-\frac{\tau}{48}\left(M_{n+1}+M_{n-1}\right)\right]^{-1}
\left[I+\frac{\tau}{48}\left(M_{n+1}+M_{n-1}\right)\right]e^{\frac{\tau}{2}Q_n}\nonumber
\end{eqnarray}
becomes a unitary matrix that conserves quadratic energy~(\ref{conser}). Since the spectrum of matrices~$Q_n$ is purely imaginary, the expression with square brackets is the Cayley formula.
The transmission matrix (\ref{schpsiMS}) was obtained using a transformation of variables, and it conserves the energy for the real spectral parameters; therefore the
corresponding scheme will be called as the fourth-order conservative transformed scheme (CT4).

Remark 1. The spectral parameter~$\xi$ is included only through the exponent exponents, as in the BO method; therefore, the use of the fast algorithm (FNFT) is difficult~\cite{Wahls2013}, but it is possible after exponential approximation~\cite{Prins2018a}.

Remark 2. An open question is how to use the free parameters $a$ and $b$ for computation with complex spectral parameters. Although the preservation of high-frequency oscillations is important, for the eigenvalues near the imaginary axis, another criterion for the scheme may be needed.

\begin{figure}[b]
\centering
\includegraphics[width=0.6\linewidth]{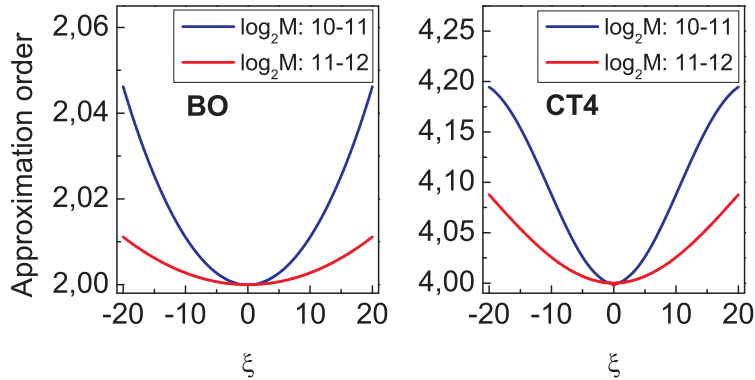}
\caption{The approximation order of the Boffetta-Osborne (BO) and conservative transformed schemes (CT4).}
\label{fig:1}
\end{figure}
The following formula was used to calculate the approximation order $m$:
\begin{equation}\label{order}
m=\log_{\frac{\tau_1}{\tau_2}}\frac{\left\Vert\tilde{\Psi}_1(L)\right\Vert_2}
{\left\Vert\tilde{\Psi}_2(L)\right\Vert_2}=
\frac{\log_2\frac{\left\Vert\tilde{\Psi}_1(L)\right\Vert_2}
{\left\Vert\tilde{\Psi}_1(L)\right\Vert_2}}{\log_2\frac{\tau_1}{\tau_2}},
\end{equation}
where $\tau_i$, $i=1$, $2$ are the steps of computational grids for two calculations with one spectral parameter $\zeta$ and $\tau_1 > \tau_2$, $\tilde{\Psi}_i(L)$ is a deviation of the calculated value $\Psi_i(L)$ from the exact analytical value $\bar{\Psi}_i(L)$ at the boundary point $t=L$. The calculations were carried out for different $p$-norms and showed close values for the approximation orders. However, for the Euclidean $2$-norm, the graphics were the smoothest.

The scheme~(\ref{schpsiMS}) was tested for different model signals, where the analytical expressions for spectral data were known. In particular, there were calculations for the oversoliton from~\cite{Satsuma1974} for a small number of discrete eigenvalues. However, to present our scheme (\ref{schpsiMS}), we chose calculations for one soliton, because this solution is smooth and not only spectral data are known for it, but also eigenfunctions~\cite{Hasegawa1995}.

Here we present numerical results for the best known potential $q(t) = \mbox{sech}(t)$. It has a single eigenvalue $\zeta_1 = 0.5i$, $b(\zeta_1) = -1$.
Since this potential is purely solitonic $b(\xi) = 0$, the continuous spectrum energy $E_c =-\frac{1}{\pi}\int\limits_{-\infty}^\infty\,\ln|a(\xi)|^2\,d\xi = 0$, while $a(\xi) = (\xi  - 0.5i)/(\xi + 0.5i)$.

\begin{figure}[t]
\centering
\includegraphics[width=0.5\linewidth]{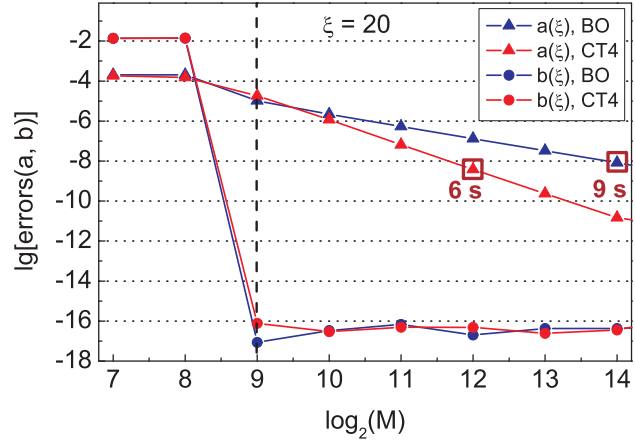}
\caption{The continuous spectrum errors.}
\label{fig:2}
\end{figure}
\begin{figure}[b]
\centering
\includegraphics[width=0.5\linewidth]{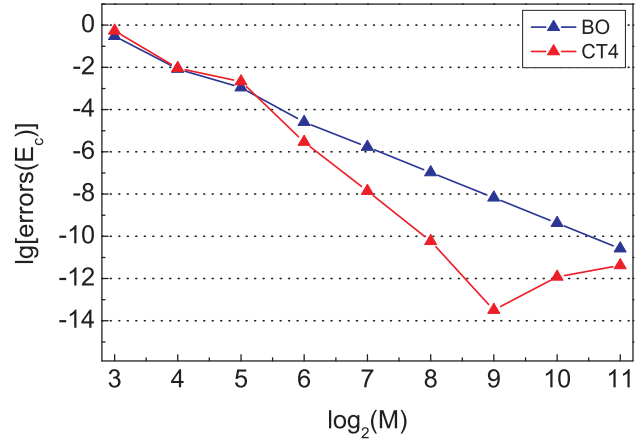}
\caption{The continuous spectrum energy errors.}
\label{fig:3}
\end{figure}

Figure~\ref{fig:1} demonstrates the approximation order~$m$ of both schemes with respect to a spectral parameter $\xi\in [-20, 20]$. Each line was calculated by the formula (\ref{order}) using two embedded grids with a doubled grid step $\tau = L/M$, $L = 40$. For the blue line, coarse and fine grids were defined by $M = 2^{10}$ and $2^{11}$. For the red one, the grid was refined one more time, namely $M = 2^{11}$ and $2^{12}$. Let us remind that the total number of points in the whole domain $[-L, L]$ is $2M +1$.

Figure~\ref{fig:2} shows the continuous spectrum errors for the fixed value of the spectral parameter $\xi = 20$.  Black dashed line in Fig.~\ref{fig:2} marks the minimum number of grid nodes~$M_{\min}$ that guarantee a good approximation. Actually, when calculating the continuous spectrum, it is necessary to choose a time step $\tau=L/M$ to describe correctly  the fastest oscillations. For a fixed value of $\xi$, the local frequency $\omega(t;\xi)=\sqrt{\xi^2+|q(t)|^2}$ of the system (\ref{psit}) varies from $\omega_{\min}=|\xi|$ to $\omega_{\max}=\sqrt{\xi^2+q_{\max}^2}$, where $q_{\max}=\max\limits_t|q(t)|$ is the maximum absolute value of the potential $q(t)$. Therefore, step $\tau$ cannot be arbitrary. In order to describe the most rapid oscillations, it is necessary to have at least 4-time steps for the oscillation period, so the inequality must be satisfied:
$$4\tau=4\frac{L}{M}\leq \frac{2\pi}{\omega_{\max}}.$$
Therefore, any difference schemes will approximate the solutions of the original continuous system (\ref{psit}) if the inequality is fulfilled for the number of points $M\geq M_{\min}=2\,L\,\omega_{\max}/\pi$.

The calculation errors for the continuous spectrum energy are compared in Fig.~\ref{fig:3}.
It is important to define the size of the spectral domain $L_{\xi}$ and the corresponding grid step $d\xi$ for the calculation of the continuous spectrum energy. According to the conventional discrete Fourier transform, we take the same number of points $N_{\xi} = N$ in the spectral domain and define a spectral step as $d\xi = \pi/(2L)$. So the size of the spectral interval is $L_{\xi} = \pi/(2\tau)$. The energy integral was computed by the trapezoid rule.

\begin{figure}[t]
\centering
\includegraphics[width=0.5\linewidth]{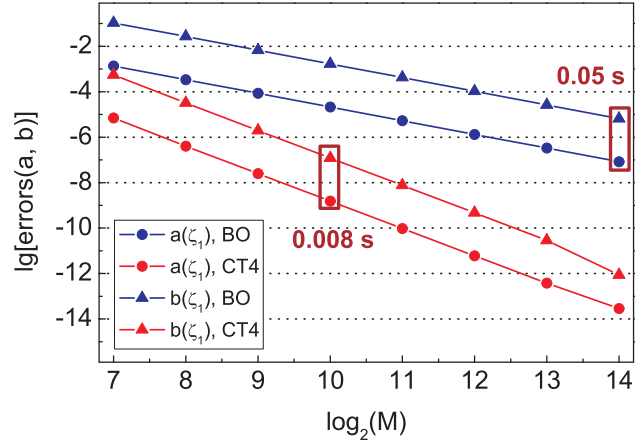}
\caption{The discrete spectrum errors.}
\label{fig:4}
\end{figure}

The discrete spectrum errors are presented in Fig.~\ref{fig:4}. The parameters $a(\zeta)$ and $b(\zeta)$ were computed for the analytically known eigenvalue $\zeta_1 = 0.5i$. In this test, we did not use any numerical algorithm to find the eigenvalue but compute  $a(\zeta)$ and $b(\zeta)$ at the exact point $\zeta = \zeta_1$ right away. It was made intentionally to estimate the error of the scheme itself and to avoid the influence of the other numerical algorithm errors.

All the errors in Figs. ~\ref{fig:2}--\ref{fig:4} are calculated using the Euclidean 2-norm.

Figures~\ref{fig:2}, \ref{fig:4} also demonstrate a comparison of the computational time. One can see that CT4 scheme allows getting a better accuracy faster than the BO scheme.

In this Letter, we proposed the family of fourth-order finite-difference one-step schemes for solving the direct Zakharov-Shabat problem on a uniform grid. Among this family, a quadratic integral preserving scheme for the continuous spectrum was distinguished. Numerical experiments for the soliton potential confirmed the theoretical order of approximation and demonstrated a significant advantage of our conservative scheme over the Boffetta-Osborne scheme. The proposed scheme works for uniform grids that can be useful when processing optical signals recorded at the receiver at regular time intervals.\\

\noindent{\bf Funding.} This work was supported by the Russian Science Foundation (grant No.~17-72-30006).

\bibliographystyle{unsrt}


\begin{thebibliography}{10}

\bibitem{ZakharovShabat1972}
V.~E. Zakharov and A.~B. Shabat.
\newblock {Exact Theory of Two-Dimensional Self-Focusing and One-Dimensional
  Self-Modulation of Waves in Non-Linear Media}.
\newblock {\em Journal of Experimental and Theoretical Physics}, 34(1):62--69,
  1972.

\bibitem{Hasegawa1973a}
Akira Hasegawa and Frederick Tappert.
\newblock {Transmission of stationary nonlinear optical pulses in dispersive
  dielectric fibers. I. Anomalous dispersion}.
\newblock {\em Applied Physics Letters}, 23(3):142--144, 1973.

\bibitem{Yousefi2014II}
Mansoor~I Yousefi and Frank~R Kschischang.
\newblock {Information Transmission Using the Nonlinear Fourier Transform, Part
  II: Numerical Methods}.
\newblock {\em IEEE Transactions on Information Theory}, 60(7):4329--4345,
  2014.

\bibitem{Turitsyn2017Optica}
Sergei~K. Turitsyn, Jaroslaw~E. Prilepsky, Son~Thai Le, Sander Wahls, Leonid~L.
  Frumin, Morteza Kamalian, and Stanislav~A. Derevyanko.
\newblock {Nonlinear Fourier transform for optical data processing and
  transmission: advances and perspectives}.
\newblock {\em Optica}, 4(3):307, 3 2017.

\bibitem{Vasylchenkova2017a}
A~Vasylchenkova, J.E. Prilepsky, D~Shepelsky, and A~Chattopadhyay.
\newblock {Direct nonlinear Fourier transform algorithms for the computation of
  solitonic spectra in focusing nonlinear Schr{\"{o}}dinger equation}.
\newblock {\em Communications in Nonlinear Science and Numerical Simulation},
  68:347--371, 3 2019.

\bibitem{Boffetta1992a}
G.~Boffetta and A.R Osborne.
\newblock {Computation of the direct scattering transform for the nonlinear
  Schroedinger equation}.
\newblock {\em Journal of Computational Physics}, 102(2):252--264, 10 1992.

\bibitem{Burtsev1998}
S~Burtsev, R~Camassa, and I~Timofeyev.
\newblock {Numerical Algorithms for the Direct Spectral Transform with
  Applications to Nonlinear Schr{\"{o}}dinger Type Systems}.
\newblock {\em Journal of Computational Physics}, 147(1):166--186, 11 1998.

\bibitem{Wahls2013}
Sander Wahls and H.~Vincent Poor.
\newblock {Introducing the fast nonlinear Fourier transform}.
\newblock In {\em 2013 IEEE International Conference on Acoustics, Speech and
  Signal Processing}, pages 5780--5784. IEEE, 5 2013.

\bibitem{Prins2018a}
Peter~J Prins and Sander Wahls.
\newblock {Higher Order Exponential Splittings for the Fast Non-Linear Fourier
  Transform of the Korteweg-De Vries Equation}.
\newblock In {\em ICASSP, IEEE International Conference on Acoustics, Speech
  and Signal Processing - Proceedings}, number~4, pages 4524--4528. IEEE, 2018.

\bibitem{Blanes2017}
Sergio Blanes, Fernando Casas, and Mechthild Thalhammer.
\newblock {High-order commutator-free quasi-Magnus exponential integrators for non-autonomous linear evolution equations}.
\newblock {\em Computer Physics Communications}, 220:243--262, 2017.

\bibitem{Prins2018}
Shrinivas Chimmalgi, Peter~J Prins, and Sander Wahls.
\newblock {Fast Nonlinear Fourier Transform Algorithms Using Higher Order
  Exponential Integrators}.
\newblock {\em arXiv preprint arXiv:1812.00703}, 12 2018.

\bibitem{Ablowitz1981}
Mark~J. Ablowitz and Harvey Segur.
\newblock {\em {Solitons and the Inverse Scattering Transform}}.
\newblock Society for Industrial and Applied Mathematics, 1981.

\bibitem{Satsuma1974}
Junkichi Satsuma and Nobuo Yajima.
\newblock {B. Initial Value Problems of One-Dimensional Self-Modulation of
  Nonlinear Waves in Dispersive Media}.
\newblock {\em Progress of Theoretical Physics Supplement}, 55:284--306, 1974.

\bibitem{Hasegawa1995}
{Hasegawa Akira} and Yuji. Kodama.
\newblock {\em {Solitons in optical communications / Akira Hasegawa and Yuji
  Kodama}}.
\newblock Clarendon Press ; Oxford University Press Oxford : Oxford ; New York,
  1995.

\end{thebibliography}

\end{document}